\documentclass[letterpaper, 10 pt, conference]{ieeeconf}
\IEEEoverridecommandlockouts                              %
\usepackage[table]{xcolor}
\usepackage[hyphens]{url}
\usepackage{graphicx}
\usepackage{epsfig}
 \usepackage{amsmath}
\usepackage{amssymb}
\usepackage{amsfonts}
\usepackage{amsthm}
 \usepackage{epstopdf}
\usepackage{algorithm}
\usepackage{algpseudocode}
\usepackage{cite}
\usepackage{subfigure}
\usepackage{afterpage}
\usepackage{float}
\usepackage{subfigure}
\usepackage{multirow}

\newcommand{\Rs}{\mathbb R}

\title{\LARGE \bf
Data-driven robust UAV position estimation\\
in GPS signal-challenged
environment
}

\author{Shenglun Yi,  Xuebo Jin,  Zhengjie Wang, Zhijun Liu and Mattia Zorzi \thanks{S. Yi  and M. Zorzi are with the Department of Information Engineering, University of Padova, Via Gradenigo 6/B, 35131 Padova, Italy. X. Jin is with the Department of Computer and Artificial Intelligence, Beijing  Technology and Business   University, No. 33 Fucheng Road, Beijing, China. Z. Wang is with the School of Mechanical and Electrical, Beijing Institute of Technology, Via No. 5, South Street, Zhongguancun, Beijing, China. Z Liu is with Beijing S-Svehicle technology.,  CO. LTD, Beijing, China.
 Emails: {\tt\small shenglun@dei.unipd.it},
 {\tt\small jinxuebo@btbu.edu.cn},
 {\tt\small wangzhengjie@bit.edu.cn}
 {\tt\small tigerzhijun@126.com},
 {\tt\small zorzimat@dei.unipd.it}}
}

\begin{document}

\maketitle
\thispagestyle{empty}
\pagestyle{empty}

%%%%%%%%%%%%%%%%%%%%%%%%%%%%%%%%%%%%%%%%%%%%%%%%%%%%%%%%%%%%%%%%%%%%%%%%%%%%%%%%
\begin{abstract}  In this paper, we consider a  position estimation problem  for an unmanned aerial vehicle (UAV) equipped with both proprioceptive sensors, i.e. IMU,  and exteroceptive sensors, i.e. GPS  and a barometer. We propose a  data-driven position estimation approach based on a robust estimator which takes into account that the UAV model is affected by uncertainties and thus it belongs to an ambiguity set.  We propose an approach to learn  this ambiguity set from the data.
\end{abstract}

%%%%%%%%%%%%%%%%%%%%%%%%%%%%%%%%%%%%%%%%%%%%%%%%%%%%%%%%%%%%%%%%%%%%%%%%%%%%%%%%
\section{Introduction}
UAV navigation systems are designed to direct a drone towards its intended destinations along paths that are both collision-free and efficient, while operating autonomously without the need for human intervention.
In the navigation task, the position estimation for the drone holds paramount importance. The latter is performed by integrating data from proprioceptive motion sensors, such as inertial measurement units (IMU) containing accelerometers and gyroscopes, with data from exteroceptive sensors, such as GPS receivers, barometer, cameras, and laser telemeters \cite{yang2021uav,7039518}.
Among them, the IMU  is not subject to interference from the external environment, while the reliability of exteroceptive sensors highly depends on the environment.

The GPS receiver, which provides users with  positioning and timing services worldwide, is the most widely used  exteroceptive sensor due to its advantages such as low cost, portability, and minimal consumption of additional computing power. Recently, there has been a surge of interest in considering the navigation task in GPS signal-challenged environments \cite{zhang2023dual,ren2017stand,tang2018vision}. This is due to the growing demand for navigation in areas where signal reception is challenging, such as indoors, under dense tree cover, or in urban canyons.
In such GPS signal-challenged environments, it has become commonplace for researchers to employ a combination of multiple sensors to enhance UAV navigation accuracy in a complementary manner. For instance, they may utilize WIFI and Ultra-Wideband (UWB) for indoor environments, or cameras for addressing challenges in both indoor and outdoor settings \cite{7402522,7081772,tang2022gnss,tang2018vision}.
However, due to cost constraints, many low-cost civilian drones cannot carry additional sensors, and their onboard computers still cannot provide sufficient support beyond basic flight plan definition and operation, such as vision-based extended computing.
For this reason, it is fundamental to design robust position estimation algorithms that rely on the limited yet reliable GPS signals to estimate the drone's position.

We address the position estimation problem using a drone equipped with both proprioceptive sensors (IMU) and exteroceptive sensors, i.e. GPS and a barometer.
This problem is a state estimation problem, which can be solved by means of the extended Kalman filter (EKF) \cite{6203568}. However, its performance highly depends on the environment. In particular, in signal-challenged environments, the accuracy of the UAV model is not so much reliable (in the sense that the model does not well describe the fact that GPS-barometer data are not so accurate). In other words, we must take into account that the model can be affected by uncertainty.
A possible way to address this issue is to consider the robust extended Kalman filter (REKF) proposed in \cite{longhini2021learning} which postulates the actual model belongs to an ambiguity set, see \cite{ROBUST_STATE_SPACE_LEVY_NIKOUKHAH_2013,
STATETAU_2017,
SING_TAC,
YI2023458,
OPTIMALITY_ZORZI,ZORZI2024105732} for a comprehensive overview. The main limitation of this approach, however, is that there are no precise rules to select this ambiguity set in practical applications.

In this paper, we propose a  data-driven position estimation algorithm based on REKF where the ambiguity set is learned from the data.   We    test the proposed algorithm using a quadcopter at Anyang’s drone test site.

% In such environments, since the performance of the exteroceptive sensors heavily relies on the external unknown environment, the GPS and barometer measurements often degrade or even intermittently disappear. To solve this problem, we  develop a data-driven position estimation algorithm, in which on the one hand an ambiguity set is defined at each time to incrementally characterize the modelling uncertainties, and on the other hand a data-driven scheme is designed to learn this ambiguity set, which is influenced by the unknown environment. More precisely, we adopt the robust extended Kalman filter  (EKF) inspired by  \cite{longhini2021learning}, as it effectively accounts for modeling uncertainties.
%This robust EKF is utilized to characterize the state estimator based on the worst-case scenario within the ambiguity set, as defined by our data-driven scheme.

The outline of the paper is as follows. The aided inertial navigation model is introduced in Section \ref{sec_2}.
Section \ref{sec_pf} introduces the  problem formulation for the position estimation in the  GPS signal-challenged environment. Then, the proposed data-driven position estimation is introduced in Section \ref{section_3}. The experimental description and the corresponding test results are presented in Section \ref{sec4} and \ref{sec_6}, respectively. Finally, Section \ref{sec_5} regards the conclusions and future research directions.

%%%%%%%%%%%%%%%%%%%%%%%%%%%%%%%%%%%%%%%%%%%%%%%%%%%%%%%%%%%%%%%%%%%%%%%%%%%%%%%%
\section{aided inertial navigation model }\label{sec_2}
The dynamic of an  unmanned aerial vehicle (UAV) system can be described by the following aided inertial navigation dynamic model \cite{martin2010generalized}:{\small \begin{equation} \label{UAV MODEL1}
\left[\begin{array}{c}
\dot q\\
\dot v\\
\dot p \\
\dot{\omega}_b\\
\dot{a}_b
\end{array}\right]=\left[\begin{array}{c}
\frac{1}{2} q * (\omega_m - \omega_b) + q *\omega_n \\
g_N +q* (a_m - a_b) *q^{-1} + q* v_n* q^{-1} \\
v \\
0\\
0
\end{array}\right]
\end{equation}}% =   \begin{bmatrix}  \rm q_1 & \rm q_2 & \rm q_3 & \rm q_4  \end{bmatrix}
where multiplication, indicated by ``$*$'', is defined by the associative law, see details in  \cite{stevens2015aircraft}; $x=\begin{bmatrix}
 q & v & p & \omega_b & a_b
\end{bmatrix}^\top $ is the state; $q =\begin{bmatrix}
 q_0 & q_1 & q_2 & q_3\end{bmatrix}^\top  \in \Rs^{4}$ is the unit quaternion representing the orientation of the body-fixed frame (b-frame) with respect to  North-East-Down frame (n-frame); $v \in  \Rs^{3} $ [m/s] is the velocity vector of the center of mass of the UAV with respect to n-frame; $p \in  \Rs^{3}$ [m] is the position  vector of the center of mass of the UAV with respect to n-frame; $\omega_b \in  \Rs^{3}$ [rad/s] is the constant vector bias on gyroscopes with respect to b-frame;  $a_b\in  \Rs^{3}$ [\text{m/s$^2$}] is the constant vector bias on  accelerometers with respect to b-frame; $\omega_m\in  \Rs^{3}$ [rad/s] is the instantaneous angular velocity vector measured by the gyroscopes; $a_m \in  \Rs^{3}$ [\text{m/s$^2$}] is the specific acceleration vector measured by the accelerometers; $w_n\in  \Rs^{3}$  is zero mean white Gaussian noise with $\mathbb E [w_{n,t} w_{n,s}]=\sigma^2_{\omega}  I_3 \delta(t-s)$ where $\delta(\cdot)$ denotes the Dirac delta function; $v_n\in  \Rs^{3}$  is zero mean white Gaussian noise with $\mathbb E [v_{n,t} v_{n,s}]=\sigma^2_{v}  I_3 \delta(t-s)$; $g_N  \in  \Rs^{3}$ [\text{m/s$^2$}]  is the constant gravity vector in the n-frame. Note that, the dynamics of $q$ and $v$ in (\ref{UAV MODEL1}) can be equivalently written in terms of the nonlinear rotational kinematics as (c.f. \cite{6203568})
 $$\begin{aligned}  \dot v &=  \text{R} (a_m - a_b+v_n) + g_N, \\
    \dot{\text{R}} &=  \text{R} S(0.5(\omega_m - \omega_b)+w_n)
 \end{aligned}$$
where $ \text{R} \in \Rs^{3\times 3}$ is the rotation matrix from  b-frame to n-frame, see details in \cite{murray2017mathematical}; $ S(x) \in \Rs^{3\times 3}$ is a skew-symmetric matrix such that $S(x)y = x \times y$  where $x,y$ are vectors, i.e. the  cross product.

\section{Problem formulation} \label{sec_pf}
Consider the problem to estimate the position of the UAV, i.e. $p$, in real-time.
We assume that the sensors available in the UAV system are: 1) An inertial measurement unit (IMU) onboard the UAV, which comprises two types of triaxial sensors that provide vector measurements expressed in the b-frame: an accelerometer that measures specific force $a_m$, and a rate gyro that measures angular velocity $w_m$; 2) A GPS receiver that gathers the first two components of the position measurements, referred to as ${ p}_{m,N}$ and ${p}_{m,E}$, expressed in the n-frame as well as the corresponding velocities; 3) A barometer that measures the third component of the position measurements, referred to as ${p}_{m,D}$, also expressed in the n-frame as well as the corresponding velocity. Therefore, the measurement model is expressed by:
\begin{equation} \label{UAV MODEL2}
y_k = Cx_k+\epsilon_k
\end{equation}
where
$$C =\left[ \begin{array}{ccccc}
0 & I_3 &  0 & 0 & 0\\
0 & 0 & I_3 & 0 & 0\\
\end{array}\right]$$
and $\epsilon_k$ is zero mean white Gaussian noise with variance  $R\in\mathbb R^6$. The latter describes the fact that the GPS measurements do not provide the exact positions and velocities. The accuracy of such measurements is described by $R$. Finally, discretizing (\ref{UAV MODEL1}) with sampling time \(\Delta\), which is the same used for the measurements, we obtain the following state-space model which describes the UAV system (including the GPS-barometer measurements):

 \begin{equation}\label{nomi_mod}
     \begin{cases} x_{k+1}= f(x_k,u_k,\varepsilon_k) + \tilde \epsilon_k\\
 y_k = C x_k +  \epsilon_k\\
 \end{cases}
 \end{equation}
where $x_k = \begin{bmatrix}  q_k^\top & v_k^\top & p_k^\top & \Delta_{\omega,k}^\top & \Delta_{a,k}^\top \end{bmatrix}^\top \in \mathbb{R}^{16}$  is the state vector, in which $q_k, v_k, p_k $ are the sampled versions of $q, v, p$, $\Delta_{\omega,k}=w_b  \Delta $ is the constant bias of the angular increment, $\Delta_{a,k} = a_b  \Delta$ is the constant bias of the velocity increment; $u_k = \begin{bmatrix}  \omega_{m,k}^\top & a_{m,k}^\top \end{bmatrix}^\top \in \mathbb{R}^{6\time2}$ is the sampled version of $[\;\omega_m^\top\; a_m^\top\;]^\top$ collected from IMU;  The noise $\varepsilon_k \in \Rs^6$ is the discretized version of $[\,w_n\;v_n\,]^\top$ and it is zero mean white Gaussian noise with covariance matrix
$$ Q_{\varepsilon,k} = \left[\begin{array}{cc}
                       \Delta \sigma^2_{\omega,k}  I_3  & 0 \\
                        0 & \Delta \sigma^2_{v,k}  I_3
                      \end{array}\right]; $$
$\tilde \epsilon_k$ is zero mean white Gaussian noise with covariance matrix
$$ Q_{\tilde \epsilon}   = \left[\begin{array}{cccc}
                      \kappa & 0 & 0 & 0\\
                        0 & \sigma^2_{q} I_3 &   0 & 0\\
                        0 & 0 &   \sigma^2_{p} I_3 & 0\\
                        0 & 0 & 0 & \kappa I_9
                      \end{array}\right]$$and it models the uncertainty arising from the discretization process. Notice that, the noise processes $\epsilon_k$, $\tilde \epsilon_k$ and $\varepsilon_k$ are independent. Moreover,    $ \sigma^2_{\omega,k} $ and  $\sigma^2_{v,k}$ are time-varying as they depend on the discretization of (\ref{UAV MODEL1}). The precise definitions of $f$, $ \sigma^2_{\omega,k} $ and  $\sigma^2_{v,k}$  as well as the  values for $\sigma^2_{\omega} $, $\sigma^2_{v} $, $ \sigma^2_{q}$ and $ \sigma^2_{p}$ can be found in the  Estimation and Control Library (ECL), \cite{PX4}.  Finally, we set $\kappa =10^{-10}$ in such a way that matrix $ Q_{\tilde \epsilon,k}$ is invertible.

The real-time estimation of the UAV system is a nonlinear state estimation problem. More precisely, at time step $k$ we can compute  the state prediction of $x_{k+1}$ given the sensor measurements by means of the EKF \cite{PX4}.  In principle, this procedure allows  to achieve a good prediction performance. However, the main limitation is that such performance can change significantly in the case the GPS signal is denied, i.e. we are in GPS signal-challenged environment. Therefore, model (\ref{nomi_mod}) does not always represent a reliable description of the underlying phenomenon.

%%%%%%%%%%%%%%%%%%%%%%%%%%%%%%%%%%%%%%%%%%%%%%%%%%%%%%%%%%%%%%%%%%%%%%%%%%%%%%%%
\section{ data-driven position estimation} \label{section_3}

We propose a data-driven robust estimation method for the UAV position estimation problem in  GPS signal-challenged environment.
Since the nominal in model (\ref{nomi_mod}) could not represent a reliable description of the environment, we assume that the actual model, i.e. the one providing a reliable description belongs to the ambiguity set defined at each time step $k$ as \begin{equation}
     \mathcal{B}_{c,k} = \left\{\ \tilde{\phi}_k(z_k|x_k) \ \text{s.t.} \ \mathcal{D}(\tilde{\phi}_k,{ \phi_k}) \le c   \right\},
 \end{equation}
where $z_{k}=[\, x_{k+1}\; y_k\,]^\top$; $\phi_k(z_k|x_k)$ is the transition probability density of the  nominal model  (\ref{nomi_mod}) at time $k$, while $\tilde \phi_k  (z_k|x_k)$ denotes the actual one and it is unknown. The mismatch between the two densities $\phi_k$ and $\tilde \phi_k$ is measured using the conditional Kullback-Leibler (KL) divergence:
 \begin{equation}
 \begin{aligned}
     \mathcal{D}(\tilde \phi_k, \phi_k) =\mathbb{\tilde E}&\left[\left.\ln\left(\frac{ \phi_k}{\tilde{\phi}_{k}}\right)\right|Y_{k-1}\right]        \end{aligned}
 \end{equation}
 where $Y_{k-1}:=\{ y_0\ldots y_{k-1}\}$, $\tilde{\mathbb E}[\cdot |\cdot]$ denotes the conditional expectation operator taken with respect to $\tilde \phi_k(z_k|x_k)\tilde \phi_k(x_k|Y_{k-1})$ and $ \tilde \phi_k(x_k|Y_{k-1})$  denotes the conditional density of $x_k$ given $Y_{k-1}$ of the actual model obtained in the previous time step. Finally, $c >0$ is the \textit{tolerance} and it defines the modeling mismatch budget at any time step: the lager $c$ is, the less reliable the model (\ref{nomi_mod}) is. It is worth pointing out this type of uncertainty arises when the density describing the nominal model is estimated from data according to the maximum likelihood principle, \cite{ROBUSTNESS_HANSENSARGENT_2008}.

To predict the state $x_{k+1}$ given $Y_k$ and taking into account that the actual model belongs to $\mathcal B_{c,k}$, we consider the solution to the minimax optimization problem:
\begin{equation}\label{minmax_pb}
    \hat{x}_{k+1} = \text{arg}\min_{g_k \in \mathcal G_k} \max_{\tilde{\phi}_k \in \mathcal B_{c,k}}\mathbb{\tilde E}\left[||x_{k+1}-g_k(y_k)||^2|Y_{k-1}\right]
\end{equation}
where $ \hat{x}_{k+1}$ is the prediction of $x_{k+1}$ given $Y_k$, $\mathcal G_k$ is the set of estimators for which the objective function in (\ref{minmax_pb})  is bounded for any $\tilde \phi_k \in \mathcal B_{c,k}$. It is worth pointing out that it is difficult to characterize a closed form expression for the solution of (\ref{minmax_pb}).

 In \cite{longhini2021learning}, an approximated characterization of the solution to (\ref{minmax_pb}) has been derived. This approximation consists in linearizing the process equation in  (\ref{nomi_mod}) around $\hat {x}_{k|k}$, which represents the estimate of $x_k$ given $Y_k$:
  \begin{equation}\label{linearized}
     \begin{array}{rl}  x_{k+1} &=   A_{k}  x_k -  A_{k} \hat {x}_{k|k} +  f(\hat{x}_{k|k},u_k) + G_k \varepsilon_k + \tilde\epsilon_k
 \end{array}
 \end{equation}
where $$\begin{aligned}  A_{k} := \left. \frac{\partial f(x, u_k, \varepsilon_k)}{\partial x}\right |_{x=\hat {  x}_{k|k}}, ~~ G_{k} := \left. \frac{\partial f(x_k, u_k, \varepsilon)}{\partial \varepsilon}\right |_{\varepsilon=\varepsilon_k}.\end{aligned}$$
Then, the robust prediction $\hat x_{k+1} $ satisfies an extended Kalman-like recursion of the form:
\begin{equation}\label{rekr}
\begin{aligned}
    L_{k} &=   V_{k}  C_{k}^{\top}( C  V_{k}  C^{\top} + R)^{-1} \\
  \hat{ x}_{k|k} &= \hat{ x}_{k} +  L_{k}( y_t - C \hat{ x}_{k} )\\
  \hat{ x}_{k+1} &=   f(\hat{ x}_{k|k},u_k) \\
    P_{k+1} &= A_{k} ( V^{-1}_{k} +  C^{\top} R^{-1} C)^{-1} A_{k}^{\top} \\
    & \hspace{3cm}  + Q_{\epsilon} +G_k  Q_{\varepsilon,k} G^\top_k \\
    V_{k+1} &= ( P^{-1}_{k+1} - \theta_k I)^{-1}
\end{aligned}
\end{equation}
where $ \theta_k > 0$ is the unique solution to $\gamma( P_{k+1}, \theta_k) = c$ and
\begin{equation} \label{theta}
    \gamma(P,\theta) = \frac{1}{2}\left\{\log\det\left(I -\theta P) + \text{tr}\left[(I-\theta P)^{-1} - I\right]\right)\right\} .
\end{equation} Notice that, the computation of $\theta$ can be efficietly accomplished using a bisection method, see Algorithm 2 in \cite{RS_MPC_IET}.
The above estimator is called robust extended Kalman filter (REKF).

The selection of the tolerance $c$ characterizing the ambiguity set has a great impact on the performance of REKF. In particular, its  optimal value depends on how much reliable the nominal model in (\ref{nomi_mod}) is. In the navigation task the model reliability depends on the particular  environment in which the
UAV system operates.
Therefore, we cannot fix $c$ a priori.

Our  main  contribution is to propose an approach for learning the tolerance through a data-driven strategy.
%Since the maximum output data rate of IMU is much higher than GPS and barometer, we apply the interpolation method to align GPS with IMU data.  The corresponding aligned GPS data, named, $x_{k}$ with $k = 0,1,\cdots,N$, is regarded as our real trajectories. At the same time, we simulate the measurements in the GPS signal-challenged environment. Specifically, we assume that the GPS signal will be intermittently denied for several seconds, while the remaining reliable GPS data that can be collected is called $\tilde y_t$ with $t = 0,1,\cdots,T$, and $T$ is much smaller than $N$. Note that for alignment purposes, we assume that the missing GPS signal is equal to the last receivable signal, and the corresponding aligned GPS, named, $y_{k}$ with $k = 0,1,\cdots,N$, is regarded as our measurements.
More precisely, given the dataset $\{(y_0,u_0) \ldots (y_N,u_N)\}$ collected from the actual model and finite set  $\mathcal{C}:=\{ c_1\ldots c_l\}$ of possible value for the tolerance, then an estimate of $c$ is given by
 \begin{equation}\label{pbc}
    \hat c = \arg\min_{c \in \mathcal{C}} \,\frac{1}{N}\sum_{k=1}^{N} \| y_{k}-\hat y_{c, k} \|^2
\end{equation}
 where $\hat y_{c, k}=C\hat x_{c, k}$ is the output prediction, while $\hat x_{c, k}$ is the state prediction using REKF with tolerance $c \in \mathcal{C}$.

Therefore, we propose the following data-driven robust position estimation algorithm:
\begin{itemize}
\item \textbf{Training phase:}
\begin{itemize}
\item for $k=1\ldots N$, compute the state prediction using REKF with a tolerance given a priori; if there is no a priori information, take it equal to zero, i.e.  compute the state prediction using EKF.
\item  for $k=N$, learn the ambiguity set, i.e. find $\hat c$ solving (\ref{pbc}).
\end{itemize}
\item  \textbf{Validation phase: }for $k>N$, compute the state prediction using REKF with tolerance $\hat c$.
\end{itemize}
In order to solve (\ref{pbc}) we have to compute $\hat y_{c,k}$ for $c\in\mathcal C$. These predictions as well as $\|y_k -\hat y_{c,k}\|^2$, with $c\in \mathcal C$, can be computed  at each time step during the first phase.  In this way, at time step $k=N$  we only have to compute $\hat y_{c,N}$, $\|y_N -\hat y_{c,N}\|^2$ and select the tolerance minimizing the loss function in (\ref{pbc}).

\section{Experiment description}\label{sec4}
\begin{figure}[t]
  \centering
  \includegraphics[width=3in]{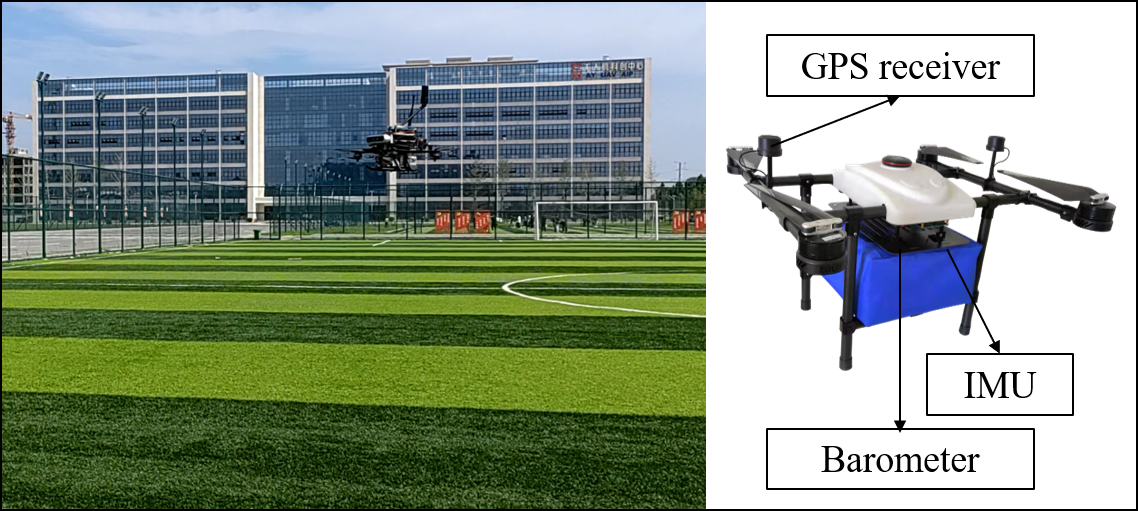}
  \caption{{\em Left.} Test site in which the flight test has been performed. {\em Right.}  The vehicle used for the flight test.}\label{UAV1}
\end{figure}

The experimental platform, as shown in the right side in Fig. \ref{UAV1}, is a quadcopter with dimensions of $35cm \times 28cm \times 25cm$ (with propeller) and mass of 1.8kg (with battery). The vehicle is equipped with Holybro Pixhawk 6C autopilot, onboard computer, IMU (with sampling frequency equal to 50Hz), GPS receiver (with sampling frequency  equal to 5Hz) and barometer (with sampling frequency  equal to 5Hz). Note that, under the ideal environment, the position accuracy of our GPS receiver  is 0.01m +  1ppm CEP. The details of the hardware configuration are summarized in Table \ref{UAV2}.

\begin{table}[h]
\centering
\caption{Hardware details.}
\begin{tabular}{|l|l|}
\hline
\rowcolor[HTML]{FFFFFF}
Component         & Specification                     \\ \hline
Platfrom          & Quadcopter, 35cm×28cm×25cm, 1.8kg \\ \hline
Autopilot         & Holybro Pixhawk 6C                \\ \hline
IMU               & ICM-42688-P                       \\ \hline
GPS               & CUAV C-RTK 9Ps                    \\ \hline
Barometer         & MS5611                            \\ \hline
\rowcolor[HTML]{FFFFFF}
On-board computer & NUC 13 Pro i5-1340P 16GB DDR4     \\ \hline
\end{tabular}
\label{UAV2}
\end{table}

\begin{figure}[t]
  \centering
  \includegraphics[width=2.0in]{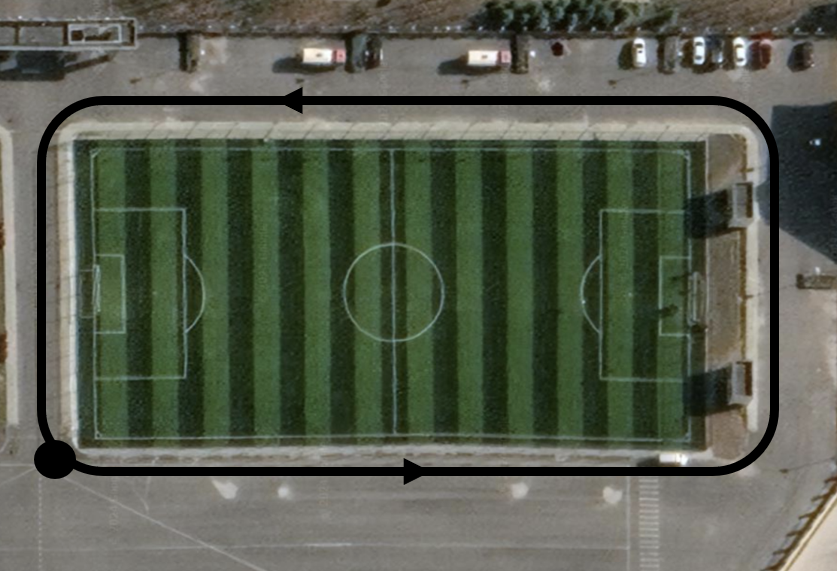}
  \caption{The vehicle started at the black circle, moving
counterclockwise, following the black path.}\label{map}
\end{figure}

\begin{figure}[t]
  \centering
  \includegraphics[width=2.3in]{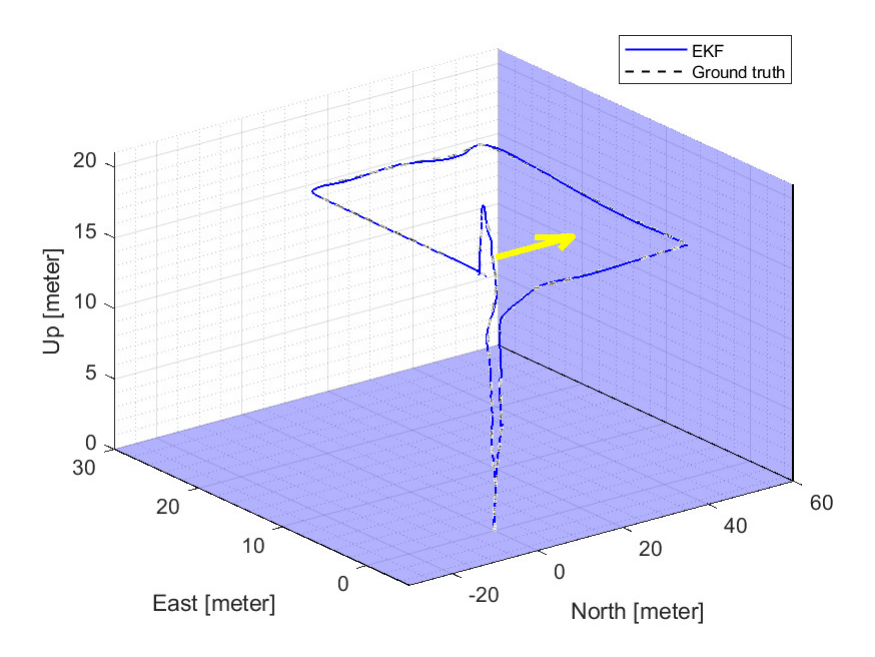}
  \caption{ Three dimensional trajectory plot for flight test  in the ideal environment: Reference trajectory (black dashed line);    Predicted position using EKF (blue line).}\label{3D1}
\end{figure}

We collected IMU, GPS, and barometer data from the flight test on an open soccer field, named Anyang's drone test site (located in Anyang, Hebei, China), under favorable weather conditions, see the left side of Fig. \ref{UAV1}. The test trajectory is illustrated in Fig. \ref{map}, with a total test duration of 106 seconds, resulting in $ T=5301$ acceleration and angular velocity measurements, as well as 531 position and velocity measurements during this test. Since the aforementioned test was conducted in an ideal environment, the position measurement error is approximately equal to 0.01 meters. Hence, we utilize the position measurements acquired in such conditions as our reference data. It is worth noting that since the maximum output data rate of the IMU is much higher than that of the GPS and barometer, we apply the linear interpolation method to align the position measurements with the IMU data. In this way, we obtain $T=5301$ reference position measurements. The latter are considered as our ground truth and are depicted in Fig. \ref{3D1} (black dashed line).

\begin{figure*}[tbp]
	\centering
	\subfigure[$ S=6s$] {\includegraphics[width=.28\textwidth]{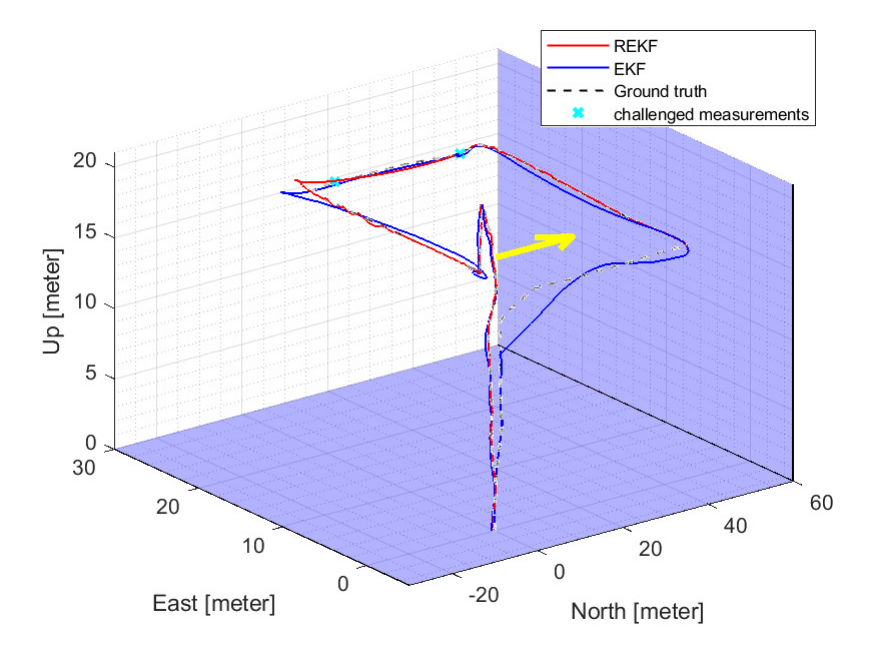}}
	\subfigure[$ S=8s$] {\includegraphics[width=.28\textwidth]{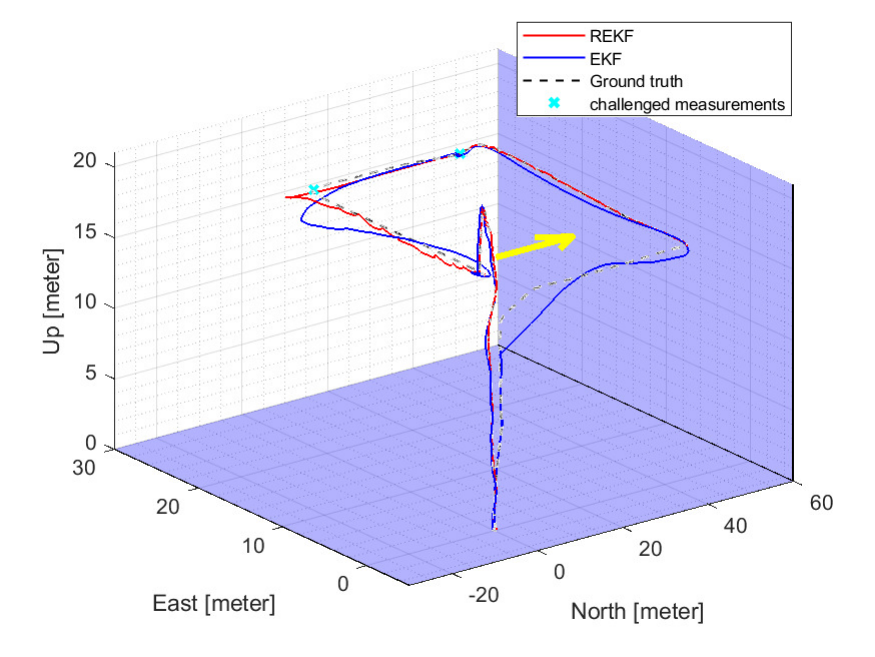}}
	\subfigure[$ S=10s$] {\includegraphics[width=.28\textwidth]{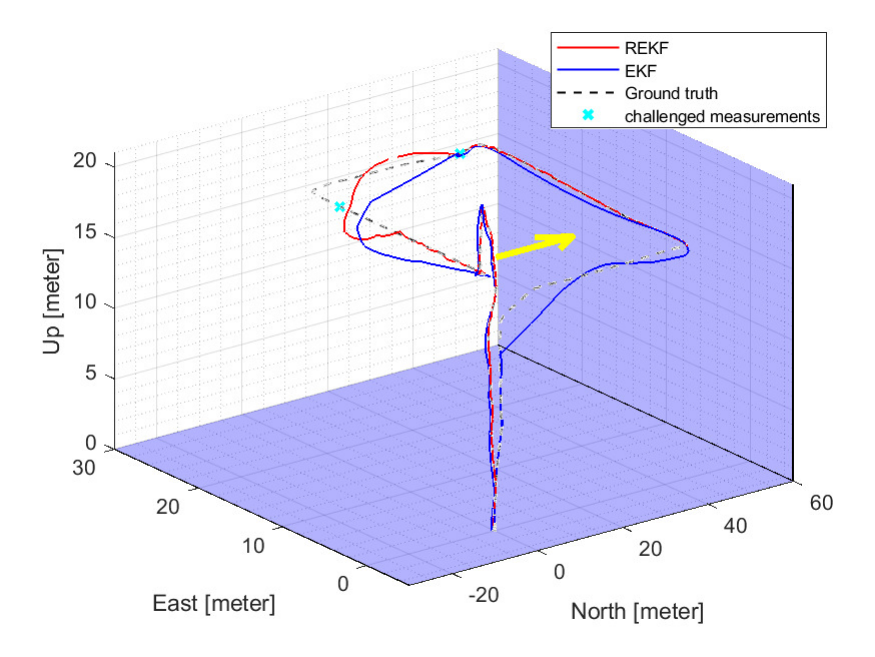}}
	\caption{ Training with $\bar S=8s$ and  validation with straight line with different GPS denial time $S$. Three dimensional trajectory plots:  Reference trajectory (black dashed line);  Predicted position using REKF (red line);  Predicted position using EKF (blue line);   GPS-barometer denied segment (between
	cyan crosses).}
	\label{fig_3D1}
\end{figure*}

\begin{figure*}[htbp]
	\centering
	\subfigure[$S=6s$] {\includegraphics[width=.28\textwidth]{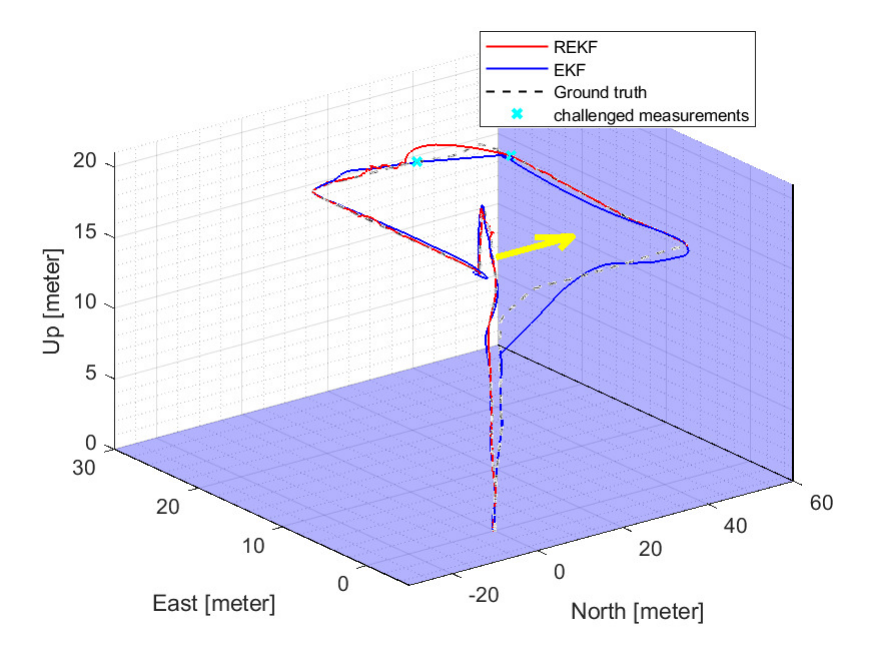}}
	\subfigure[$ S=8s$] {\includegraphics[width=.28\textwidth]{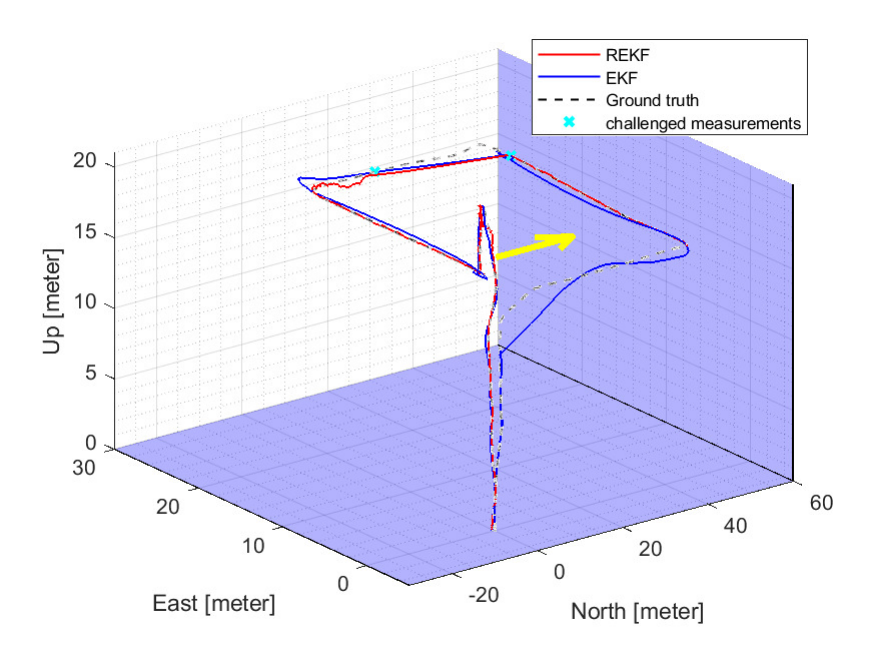}}
	\subfigure[$ S=10s$] {\includegraphics[width=.28\textwidth]{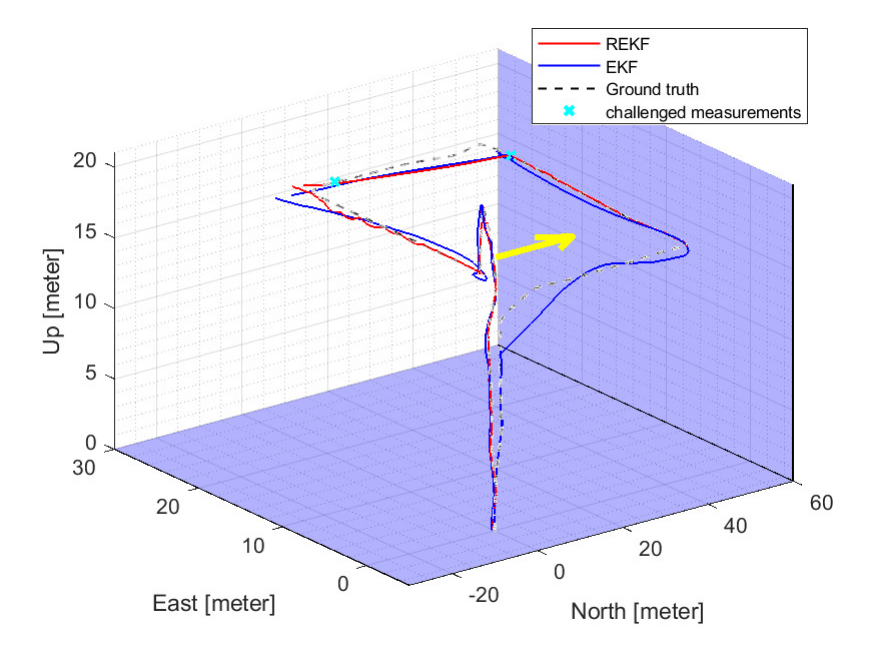}}
	\caption{Training with $\bar S=8s$ and  Validation with turn with different GPS denial time $S$.  Three dimensional  trajectory  plots: Reference trajectory (black dashed line);  Predicted position using REKF (red line);  Predicted position using EKF (blue line);   GPS-barometer denied segment  (between cyan crosses).}
	\label{fig_3D2}
\end{figure*}

From this experiment we  create two different data sets corresponding to different scenarios.\\
\textbf{GPS signal-normal environment.} We consider the case in which the GPS-barometer signal is always available. Since the maximum output data rate of the IMU is much higher than that of the GPS and barometer, we apply the causal zero-order hold interpolation to align position and velocity measurements with IMU data. The latter represents our measurements data. \\
%\textbf{\ys Intermittently GPS-denied  environment.} {\ys We consider that the flight mission in an intermittently GPS-denied  environment.} More precisely, we assume that only one reliable position signal can be collected in every 4 seconds, leaving only 27 receivable position and velocity measurements. In what follows, we refer to the corresponding position measurements as challenged measurements, which are depicted in Fig. \ref{3D2} using green crosses.
%Also, for alignment purposes, we assume that the missing GPS-barometer signal points are equal to the last available one. In this way, we obtain the corresponding 5301 challenged position and velocity measurements, named $\ys \bar Y_T = \{ \bar y_1, \cdots, \bar y_k, \cdots, \bar y_T \}$. Note that, $\ys \bar Y_T$ can be viewed as a zero-order hold interpolation obtained form the  27 receivable position and velocity measurements.\\
\textbf{GPS signal-challenged environment.}  We divide the dataset into two  segments: the first part corresponding to the training phase and the second one to the validation phase. The first segment starts from the moment the drone takes off and its duration  is equal to 30 seconds. More precisely, we simulate three different training segments scenarios:
\begin{itemize}
  \item \textbf{Training with $\mathbf{\bar S=6s}$:}  from takeoff, the denied GPS-barometer signal is from the 15th second to 21th second, i.e. lasts for $\bar S=6s$.
  \item  \textbf{Training with $\mathbf{\bar S=8s}$:} from takeoff, the denied GPS-barometer signal is from the 15th second to 23th second, i.e. lasts for $\bar S=8s$.
  \item  \textbf{Training with $\mathbf{\bar S=10s}$:} from takeoff, the denied GPS-barometer signal is from the 15th second to 25th second, i.e. lasts for $\bar S=10s$.
\end{itemize}
The second segment starts from the 30th second and its duration is equal to 76 seconds. More precisely, we consider the following two types of validation segments:
\begin{itemize}
  \item  \textbf{Validation with straight line.} The drone loses GPS signal while flying in a straight line, i.e. starting from the 40th second and lasts for $S$ seconds, see e.g. the GPS-barometer denied segment between two cyan crosses in Fig. \ref{fig_3D1};
  \item \textbf{Validation with turn.} The drone loses GPS signal just as it is about to make a turn, starting from the 36th second and lasts for $S$ seconds,  see e.g. the GPS-barometer denied segment between two cyan crosses in Fig. \ref{fig_3D2}.
\end{itemize}
Note that, during the GPS-barometer denied segment, we force all missing GPS-barometer signal points to remain constant and equal to the last available value.
Finally, we apply the causal zero-order  hold interpolation to align the corresponding position and velocity measurements with IMU data. The latter represents our measurements data.

\section{TEST RESULTS}\label{sec_6}
In this section, we want to test the performance of the data-driven robust position estimator proposed in Section \ref{section_3}, hereafter referred to as REKF. More precisely, we compare it with the extended Kalman filter, hereafter referred to as EKF.

%\begin{figure*}[t]
%  \centering
%  \includegraphics[width=7.3in]{10s_sub-eps-converted-to.pdf}
%  \caption{Trajectory along the three directions in the second phase. From left to right:  Position on North, East and Up: Ground truth (black dashed line);  Predicted position using EKF (blue line); Predicted position using REKF (red line).  }\label{1D}
%\end{figure*}

The state prediction is performed according  model (\ref{nomi_mod}) whose sampling time is $\Delta=0.02$ seconds (i.e. the one of the IMU signal).
Both REKF and EKF are initialized with
$$ \begin{aligned}\hat x_0 =  &[ 1, ~ 0, ~ 0, ~ 0,  ~0, ~ 0, ~ 0, ~ 0.9926, ~ -0.0126, ~ 0.0230, \\
& 2.7556 \times 10^{-6}, ~ 2.7556\times 10^{-6}, ~ 2.7556\times 10^{-6}, \\
  & 6.7600\times 10^{-11}, ~ 6.7600\times 10^{-11}, ~ 6.7600\times 10^{-11}  ]^\top, \end{aligned} $$
which represent the initial prediction.
Note that, the initial quaternion is set according to the initial attitude of UAV; the initial velocity is $[ 0, ~ 0, ~ 0]$; the initial position, measured by the GPS receiver, is $[0.9926, ~ -0.0126, ~ 0.0230]$, and the initial constant vector biases are determined according to the accuracy of the accelerometer and the gyroscope, respectively, see details in \cite{PX4}. Moreover, since we know $\hat x_0$ represents a relatively good estimate of $x_0$, the initial covariance matrix of the prediction error, i.e. $V_0$, should be taken not too large. Thus, we set:
\begin{equation}V_0 =  10^{-3} \cdot I_{16}. \label{V0}\end{equation}
Then, the parameters for process noises, i.e. $Q_{\varepsilon,k}$ and $ Q_{\tilde \epsilon}$, are computed as outlined in Section \ref{sec_pf}.  The covariance matrix of the measurement noise process $\epsilon_k$, i.e. $R\in \Rs^{6\times 6}$, is assumed to be diagonal, say $R= {\rm diag}(r)$.  Since the accuracy for GPS and barometer under ideal conditions are equal to $0.7m$ and $0.8m$, respectively, the elements in the main diagonal of $R$, i.e. $r$, are  selected   in such a way the model respects the  measurement accuracy  with the $95\%$ confidence level:
 $$r= \begin{bmatrix}
 1.96 & 1.96 & 2.56 & 0.1225 & 0.1225 & 0.16\end{bmatrix}. $$
 We test EKF using the GPS signal-normal environment dataset with the
aforementioned parameters setting: as we can see in Fig. \ref{3D1} EKF achieve the very good prediction performance.
 Accordingly, in what follows we still continue to use this parameters setting.

\begin{table*}[t]
\centering
\caption{$\overline{\text{RMSE}}$ of  the prediction error using different filters and in the different cases.}
{\scriptsize
\begin{tabular}{|c|c|ll|ll|ll|}
\hline
\multicolumn{1}{|l|}{}                                                                       & \multicolumn{1}{l|}{} & \multicolumn{2}{c|}{\begin{tabular}[c]{@{}c@{}}{Training  with}\\ $\bar S=6s$\end{tabular}} & \multicolumn{2}{c|}{\begin{tabular}[c]{@{}c@{}}{Training with}\\ $\bar S=8s$\end{tabular}} & \multicolumn{2}{c|}{\begin{tabular}[c]{@{}c@{}}{Training with}\\ $\bar S=10s$\end{tabular}} \\ \hline
\multicolumn{1}{|l|}{}                                                                       & \multicolumn{1}{l|}{} & \multicolumn{1}{c|}{EKF}                   & \multicolumn{1}{c|}{REKF}                & \multicolumn{1}{c|}{EKF}                   & \multicolumn{1}{c|}{REKF}                & \multicolumn{1}{c|}{EKF}                    & \multicolumn{1}{c|}{REKF}                \\ \hline
\multirow{3}{*}{\begin{tabular}[c]{@{}c@{}}{Validation}\\ (Straight line)\end{tabular}} & $ S=6s$                  & \multicolumn{1}{l|}{0.9010}                & 0.5284                                   & \multicolumn{1}{l|}{0.9395}                & 0.5380                                   & \multicolumn{1}{l|}{0.9709}                 & 0.4465                                   \\ \cline{2-8}
                                                                                             &$   S=8s  $                & \multicolumn{1}{l|}{0.9912}                & 0.6383                                   & \multicolumn{1}{l|}{1.0194}                & 0.6744                                   & \multicolumn{1}{l|}{1.0573}                 & 0.5722                                   \\ \cline{2-8}
                                                                                             & $ S=10s  $              & \multicolumn{1}{l|}{1.1038}                & 0.7158                                   & \multicolumn{1}{l|}{1.1451}                & 0.6096                                   & \multicolumn{1}{l|}{1.1239}                 & 0.5505                                   \\ \hline
\multirow{3}{*}{\begin{tabular}[c]{@{}c@{}}{Validation}\\ (Make a turn)\end{tabular}}   & $ S=6s  $                & \multicolumn{1}{l|}{0.5991}                & 0.3748                                   & \multicolumn{1}{l|}{0.6291}                & 0.3580                                   & \multicolumn{1}{l|}{0.6514}                 & 0.3343                                   \\
\cline{2-8}
&$ S=8s $                & \multicolumn{1}{l|}{0.8306}                & 0.4877                                   & \multicolumn{1}{l|}{0.8699}                & 0.4388                                   & \multicolumn{1}{l|}{0.9013}                 & 0.3893                                   \\
\cline{2-8}
&$ S=10s   $              & \multicolumn{1}{l|}{1.1069}                & 0.6407                                   & \multicolumn{1}{l|}{1.1554}                & 0.6858                                   & \multicolumn{1}{l|}{1.1666}                 & 0.5413                                   \\ \hline
\end{tabular}
}
\end{table*}

Next, we compare REKF and EKF using the GPS signal-challenged environment dataset.   The set $\mathcal C$ contains $l=40$ elements equispaced in the interval $[2 \cdot 10^{-4},\;1]$.
In this setting, we obtain that the optimal tolerances are $\hat c = 0.0274$, $\hat c = 0.1013$ and $\hat c = 0.5199$ for  the training segments with $\bar S=6s$, $\bar S=8s$ and $\bar S=10s$, respectively. Notably, as the duration of the GPS-barometer signal denial increases, thereby increasing the uncertainty, $\hat c$ also increases.
Fig. \ref{fig_3D1} shows the 3-dimensional predicted positions corresponding to  the training segment with $\bar S = 8s$ and the validation segment with straight line for different values of $S$.  During the training phase, the prediction of REKF matches that of EKF because $c$ is set to zero.  In the validation phase,  EKF provides an acceptable performance when $S=6$s (i.e. when the GPS-barometer signal is denied for a short period). However, as $S$ increases, EKF performs poorly, even after the GPS signal is restored, particularly during UAV maneuvers such as turning and landing. In contrast, REKF performs significantly better, especially when the denied segment in the training phase coincides that in validation one, i.e. $\bar S=S=8s$.  Fig.  \ref{fig_3D2} presents the 3-dimensional predicted positions corresponding to  the  training segment  with $\bar S = 8s$ and the  validation segment with turn for different values of $S$. Such scenario  is more challenging than the first one. Indeed, during the period when the GPS-barometer signal is denied and the drone is turning,  the predicted trajectories of both EKF and REKF show  a remarkable deviation in respect to the ground truth.

Another noteworthy point of discussion is whether the proposed REKF can also improve the prediction performance when the GPS-barometer data is restored.  This aspect  has a crucial impact on the ability of UAV to escape from a GPS-barometer signal-challenge environment. As we can see in Fig. \ref{fig_3D1}-\ref{fig_3D2},  REKF  aligns in a fast way with the actual trajectory once the GPS-barometer signal is reacquired.   Next, we evaluate the  performance of EKF and REKF in the interval starting when the GPS-barometer signal is reacquired (the corresponding point is denoted by $K$) and whose ending coincides with the one of the dataset (i.e. the ending point is $T=5301$).
More precisely, we consider,  the root mean square error (RMSE) of the prediction error along  the three dimensions:
{\small \begin{equation*}
\begin{aligned}
 &\overline{\text{RMSE}} =   \sqrt{ \frac{\sum^T_{t=K+1} \| {{p}}^N_t-\hat{p}^N_{t} \|^2}{T-K}}   \\
  & \hspace{0.5cm}  + \sqrt{ \frac{\sum^T_{t=K+1} \| {{p}}^E_t-\hat{p}^E_{t} \|^2}{T-K} } + \sqrt{ \frac{\sum^T_{t=K+1} \| {{p}}^D_t-\hat{p}^D_{t} \|^2}{T-K} }
 \end{aligned}
\end{equation*}}
where   $p_t=[p_t^N,~p_t^E,~p_t^D]$ is the reference position at time $t$ extracted from ground truth, while $\hat{p}^k_{t} = [\hat p_t^N,~\hat p_t^E,~\hat p_t^D]$ is the position prediction extracted from the estimator $\hat x_t$. As we can see,   REKF  outperforms EKF, and the larger the denial time during the training phase is, the better the performance of REKF is. In this respect, we can preliminarily conclude that  REKF is able to rapidly realign the predicted trajectory with the ground truth when the GPS-barometer signal is reacquired.

\section{Conclusions}\label{sec_5}
We have faced the problem to estimate (or more precisely, predict) the position of UAV systems by means of IMU and GPS-barometer data. The latter is fundamental in UAV navigation tasks. The estimation  performance highly depends on the environment, i.e. whether GPS-barometer data are available or they are not always available. In many real situations the environment is not known a priori, a fact that can lead to a remarkable accumulated drift in the position estimation. To take into account such issue, we have proposed a data driven robust EKF estimator where in the first phase the model uncertainty (which is related to the fact that the environment is not known) is inferred from the data, then in the second phase the estimation is performed taking into account such model uncertainty. We have tested the proposed algorithm to a quadcopter showing its superiority with respect to the standard EKF algorithm. An open question we will address in the future regards whether our approach is better than adaptive EKF filters. The main weakness of the latter is that the adaptation involves the optimization of many parameters, while in ours we have to optimize only one parameter, i.e. the tolerance characterizing the ambiguity set.

%%%%%%%%%%%%%%%%%%%%%%%%%%%%%%%%%%%%%%%%%%%%%%%%%%%%%%%%%%

\bibliography{biblio}
\bibliographystyle{ieeetr}

\end{document}